\documentclass[11pt]{amsart}
\usepackage[letterpaper, portrait, margin=1in, headheight=0pt]{geometry}
\usepackage{enumerate}
\usepackage{amsmath,amsthm,amssymb,amsfonts}
\usepackage{tikz-cd}
\usepackage{nicefrac}
\usepackage{hyperref}
\usepackage{amsrefs}

\usepackage{amssymb}
\usepackage{amsthm} 
\usepackage{caption}
\usepackage{wrapfig}

\newcommand {\mm}[1] {\ifmmode{#1}\else{\mbox{\(#1\)}}\fi}

\usepackage{mathtools}

\theoremstyle{definition}
\newtheorem{theorem}{Theorem}[section]
\theoremstyle{definition}
\newtheorem{definition}[theorem]{Definition}

\theoremstyle{remark}
\newtheorem{remark}[theorem]{Remark}
\numberwithin{equation}{section}

\usepackage{cite}

\DeclareMathOperator{\mir}{mir}
\DeclareMathOperator{\Span}{Span}

\usepackage{mathtools}
\DeclarePairedDelimiter{\norm}{\lVert}{\rVert}

\renewcommand{\phi}{\varphi}

\newcommand{\CF}[2]{\mathcal{C}^{#1}_{\mathcal{F}_{\lceil #2\rceil}}}
\newcommand{\CFn}[2]{\mathcal{C}^{#1}_{\mathcal{F}_{\lceil #2\rceil_{n}}}}
\newcommand{\CFa}[2]{\mathcal{C}^{#1}_{\mathcal{F}_{\lceil #2\rceil_{a}}}}
\newcommand{\CFp}[2]{\mathcal{C}^{#1}_{\mathcal{F}_{\lceil #2\rceil_{P}}}}
\newcommand{\CFt}[2]{\mathcal{C}^{#1}_{\mathcal{F}_{\lceil #2\rceil_{T}}}}

\usepackage{color}

\begin{document}

\title[Structure of the Jones Polynomial]{Big Data Approaches to Knot Theory: Understanding the Structure of the Jones Polynomial}

\author{Jesse S F Levitt}
\address{University of Southern California, 
Los Angeles, CA USA}
\email{jslevitt@usc.edu}

\author{Mustafa Hajij}
\address{KLA corporation,
Ann Arbor, MI USA}
\email{mustafahajij@gmail.com}

\author{Radmila Sazdanovic}
\address{North Carolina State University,
Raleigh, NC USA}
\email{rsazdanovic@math.ncsu.edu}

\begin{abstract}
We examine the structure and dimensionality of the Jones polynomial using manifold learning techniques.
Our data set consists of more than 10 million knots up to 17 crossings and two other special families up to 2001 crossings.
We introduce and describe a method for using filtrations to analyze infinite data sets where representative sampling is impossible or impractical, an essential requirement for working with knots and the data from knot invariants.
In particular, this method provides a new approach for analyzing knot invariants using Principal Component Analysis.
Using this approach on the Jones polynomial data we find that it can be viewed as an approximately 3 dimensional manifold, that this description is surprisingly stable with respect to the filtration by the crossing number, and that the results suggest further structures to be examined and understood.
\end{abstract}

\maketitle

\section{Introduction}

Throughout the development of low-dimensional topology, there has been an emphasis on the study of invariants from algebraic, combinatorial and geometric perspectives.
The scarcity of results considering the statistical nature of these invariants is quite surprising given the abundance of available data. Examining the distributions that arise from invariants should reveal and illuminate structures that are difficult to see using traditional tools.
We consider the Jones polynomial from a statistical perspective and provide an outline of how to use filtrations to investigate infinite data sets of this type.

The techniques of big data and deep learning provide useful tools for analyzing the statistical nature of knot invariants.
Multiple advances have resulted from melding traditional methods in Physics and Mathematics with the emerging data-driven techniques of scientific computing.
These have ranged from solving previously intractable problems in Computer Vision to providing significant improvements in earthquake prediction models.
Despite the wide number of techniques available, the study of how to use these powerful statistical tools in pure mathematics is in its infancy.
Machine learning and data mining techniques have just started to attract attention in Knot theory, see \cite{ward2018using,jejjala2019deep,hughes2016neural} and to our knowledge have largely been used for predicting valuations.
This is the first in a series of papers examining how to apply these techniques to low-dimensional topology to gather structural insights.
We start by focusing on dimensionality reduction, with further analysis using supervised machine learning techniques~\cite{hastie2005elements} and persistent homology~\cite{edelsbrunner2000topological} forthcoming.

Low-dimensional topology, and knot theory in particular, is among the most data-rich of mathematical sub-fields.
Tabulating data concerning knots is a longstanding tradition dating back to the 1860s~\cite{Thomson1869}.
As computing power improves and people continue to search for answers to fundamental questions in knot theory such as the Jones unknot conjecture~\cite{jones1997polynomial} and the hyperbolic volume conjecture~\cite{murakami2011introduction}, people have continued to enlarge our tabulations of known knots.
Recently, Burton tabulated all the prime knots up to and including 19 crossings~\cite{BBurton18}, finding over 350 million total prime knots.
This tabulation was summarized as part of the software package Regina~\cite{regina} with published DT-codes as defined in~\cite{DTCodeIntro}.
A separate effort at tabulating large numbers of unique knots was also recently undertaken by Tuzun and Sikora in their demonstration that no counterexamples to the Jones unknot conjecture exist up to 23 crossings~\cite{SikoraTuzun18}. In their tabulation, well over 10 trillion knot diagrams were considered using distinct methods from Burton's.
While the exact number of distinct knots with a certain number of crossings is still unknown, we do know that this number grows at an exponential rate as we increase the number of crossings~\cite{ernst1987growth}.
This ensures that data-related questions arising from Knot Theory naturally fit into a big data framework. 

The first major contribution of this paper is to demonstrate a reliable technique by which manifold learning can be applied to infinite data sets where representative sampling is impossible or impractical.
We describe how to analyze and construct a usable filtration on the infinite set of knots, where the Jones polynomial is unbounded in degree.
The experimental results of this demonstrates that our Jones invariant data is well approximated by a three dimensional manifold, consistent across our filtration up to computation limits.

Meaningfully applying dimensionality reduction techniques to our data proved difficult.
It is unknown how to create a representative subsample of Jones polynomial data.
Any conclusions drawn must remain consistent when choosing comparable subsets of the data.
Results were sensitive to the choice of how encode the data for comparison.
Using the same approach as was used to find patterns in the more general coloured Jones polynomial~\cite{lee2018trivial,beirne2017q,lovejoy2013bailey,bataineh2016colored,armond2011rogers,elhamdadi2017foundations} provided data where any structures proved transient.
Fortunately, exactly one model for encoding the data provided results that were both remarkable and persistent, it is discussed in Section~\ref{sec:data}.

The requirement that results be persistent across comparable subsets of the data required detailed analysis of how to filter sets of knots into related families.
Knot Theory has always driven researchers to calculate knot invariants and organize them into data tables in a process called knot tabulation \cite{hoste1998first,hoste2005enumeration}.
Originally envisioned as a way to distinguish different atomic properties~\cite{Thomson1869}, modern work has suggested that a classification system could assist in the understanding of glueball particles~\cite{FlamminiStasiak2007}.
Since then a series of systems have been suggested for ordering or relating knots within these ever expanding tabulations~\cite{DiaoErnstStasiak2009, CantarellaHenrichEtAl2017}.

Upon generating our Jones polynomial data for all knots up to 17 crossings we examined several methods for organizing the data.
We considered the crossing number, Rasmussen s-invariant~\cite{rasmussen2010khovanov}, signature, unknotting number, and a wide variety of properties intrinsic to the Jones data itself.
As we discuss in Section~\ref{sec:results} below, organizing the data by crossing number yields persistent results despite the manner in which the set varies dramatically in both the ratio of alternating to nonalternating knots present in the sample and the expanding size of the data considered.

\begin{figure}[h]
  \centering
   {\includegraphics[width=\textwidth]{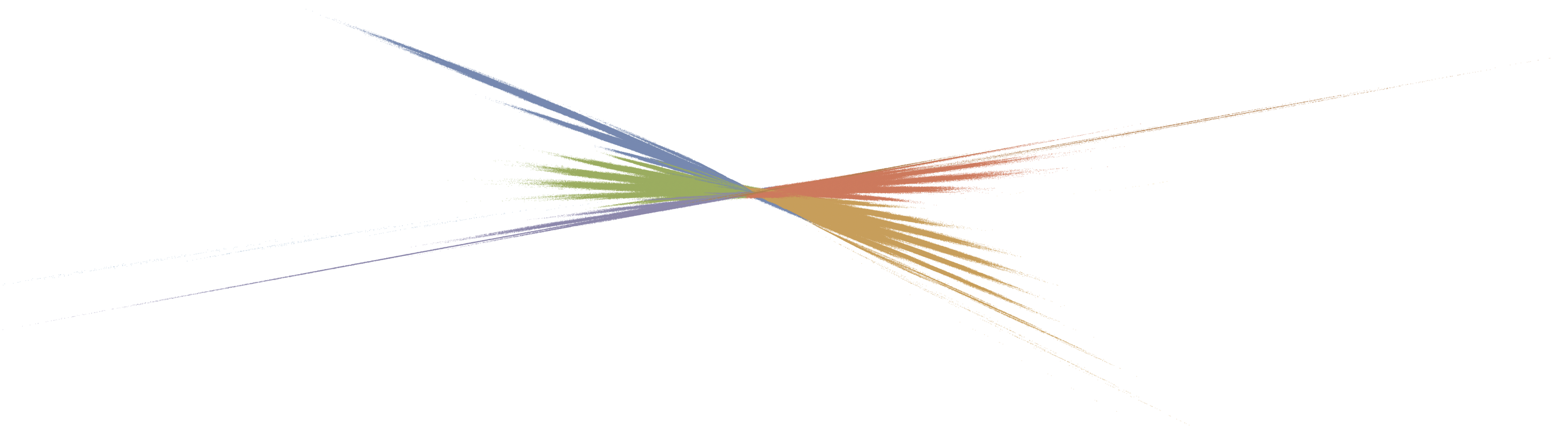}
    \caption{A PCA projection of the Jones polynomial data into 3 dimensions colored by the knot signature to highlight internal structure.}
  \label{fig:pcaJones}}
\end{figure}

The second major result from this study is to demonstrate a new tool for comparing knot invariants and understanding their structure.
Applying dimensionality reduction to the Jones data using Principal Component Analysis (PCA)~\cite{wold1987principal} as in Figure~\ref{fig:pcaJones} we see a rich three dimensional structure with large scale features differentiated by their signature with subfeatures of smaller `tendrils' with as yet unknown significance.
We propose the following definition for understanding the results of dimensionality reduction via PCA, given the discussion of Remark~\ref{rem:ha}.
\begin{definition}[Dimension of a polynomial knot invariant, $P$]
Let $\mathbf{k}$ be the value for which the normalized explained variance of the first $\mathbf{k}$ PCA components sums to more than 95\%. If this remains stable across the crossing filtration, then the knot invariant $P$ has \emph{dimension} $\mathbf{k}$.
\end{definition}
Under this definition the Jones polynomial is 3 dimensional, the $Z_0$ polynomial of Bar-Natan and van der Veen~\cite{bar2017polynomial} is 2 dimensional and the Alexander polynomial~\cite{alexander1928topological} is 1 dimensional.
In this paper we will focus solely on understanding how the Jones polynomial is 3 dimensional under this definition, while the remaining calculations will be published in upcoming work.

To perform this analysis we relied on the tools from the KnotTheory package~\cite{KA} for calculating the Jones polynomial for all knots.
The DT codes we used for knots up to 16 crossings were exactly those in the KnotTheory package, while to calculate data for the 17 crossing knots we added the DT codes from Burton's Regina program~\cite{BBurton18} to the KnotTheory data tables.
The PCA calculations were done using the scikit-learn library~\cite{scikit-learn} and in Mathematica~\cite{Mathematica}. Finally, knot figures were generated using Inkscape~\cite{Inkscape}.

\section{Background}
In this section we briefly provide the reader with an overview of the definitions and notions used in the article.
We begin with the definition of the Jones polynomial, followed by an overview of the basic properties of the PCA technique. 

\subsection{The Jones Polynomial}
\label{ssec:JP}

The Jones polynomial~\cite{jones1997polynomial} and its generalizations~\cite{murakami2001colored,reshetikhin1990ribbon,turaev1988yang} play a fundamental role in low-dimensional topology~\cite{murakami2001colored,kashaev1997hyperbolic,murakami2011introduction}.
The discovery of the Jones polynomial has  led to multiple major discoveries in various areas of low-dimensional topology \cite{le1905aj,khovanov2005categorifications,dasbach2006head,bar1996melvin,le2006colored,garoufalidis2005colored,le2000integrality}.
Understanding the discriminative power of the Jones polynomial, its relations to other classical invariants of knots and links, as well as the information encoded in its coefficients, conjectured to be related to the hyperbolic volume of the knot, are important problems in low-dimensional topology. 
Furthermore, the coefficients of the Jones polynomial and its generalizations have been proven to be related to many interesting areas in number theory, and have been the subject of an extensive research effort~\cite{lee2018trivial,elhamdadi2017pretzel,hajij2016tail,beirne2017q,lovejoy2013bailey,hajij2017colored,bataineh2016colored,armond2011rogers,elhamdadi2017foundations}.

\begin{figure}[!ht]
  \centering
  { \includegraphics[height=0.6in]{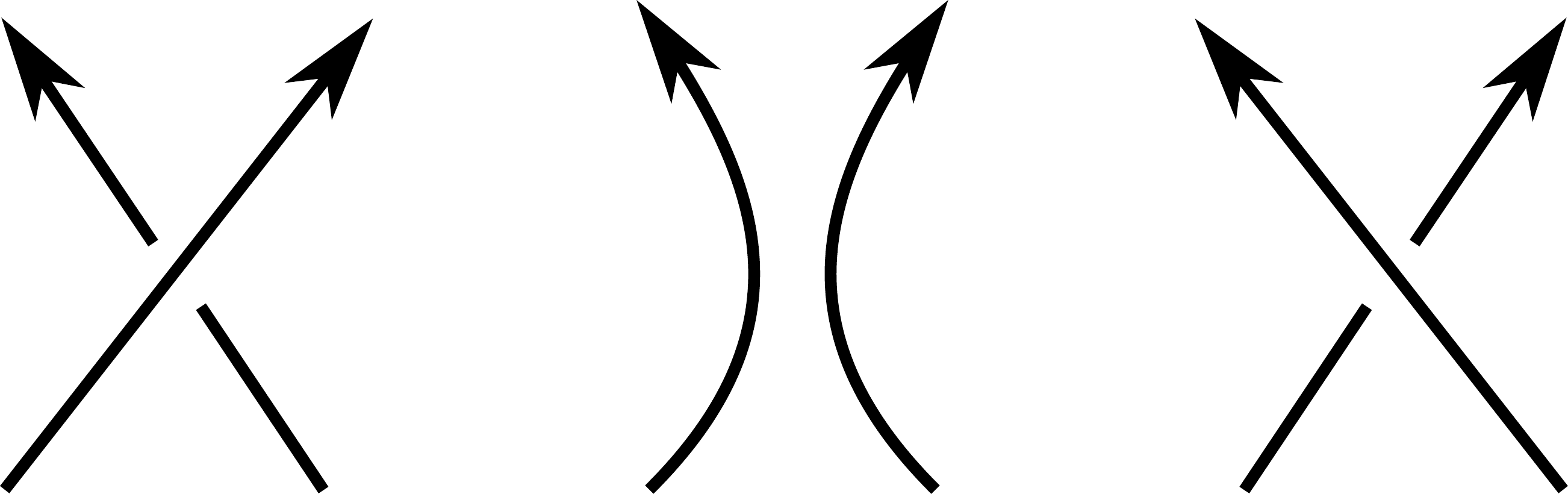} \hspace{2cm}\includegraphics[height=0.65in]{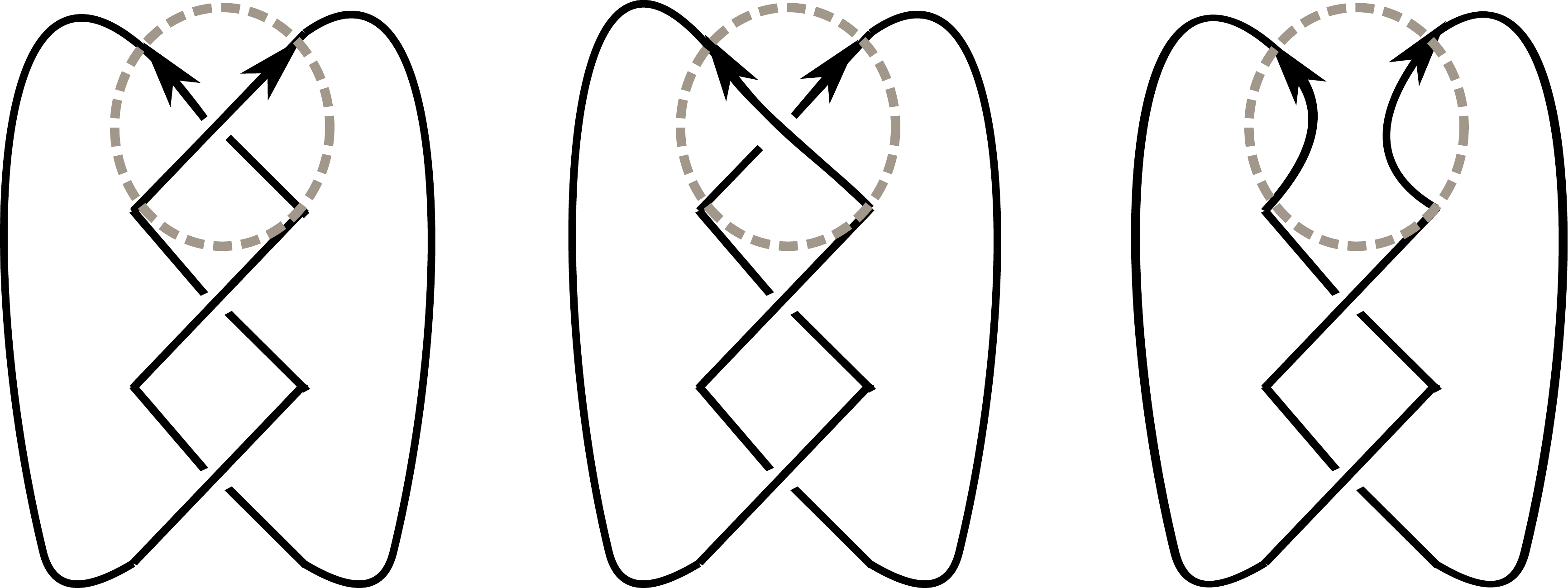}
    \caption{The three pictures on the left represent a positive crossing, a smoothing, and a negative crossing, denoted $ L_{+}$, $L_{0}$ and $L_{-}$ respectively. The right picture represents the same concept, fixed on the top crossing in the trefoil diagram $L$ first with $L_{+}$, a diagram of an unknot with $L_{-}$ next and of the Hopf link with $L_{0}$ at the end.}
  \label{3crossings}}
\end{figure}

Let $K$ be a knot in $\mathbb{S}^3$. The Jones polynomial, denoted by $J_{K}(q)$, is a Laurent polynomial in $\mathbb{Z}[q^{\pm 1}]$. The Jones polynomial can be characterized by the requirements that $J_{K}(q)=1$, when $K$ is the unknot, and that it satisfies the following skein relation:
$
(q^{1/2}-q^{-1/2})J_{L_{0}}(q)=q^{-1}J_{L_{+}}(q)-qJ_{L_{-}}(q).$
Here $L_0$, $L_{-}$ and $L_{+}$ are three oriented link diagrams that are identical everywhere except at single crossing as appears on the right in Figure \ref{3crossings}.
The skein relation can be used to compute the Jones polynomial for any given link $L$.

\subsection{Principal Component Analysis}
\label{ssec:PCA}

Principal Component Analysis (PCA) is one of the most popular multivariate statistical techniques in big data.
PCA is defined as an orthogonal linear transformation that maps a given data set $X$ to a new orthonormal basis that is aligned with the core properties of the data.
To accomplish this, PCA finds and ranks the linear directions along which the data has maximal variance.

Given the data set $X$, where $|X|=n,$ and $X\subset\mathbb{R}^d$, the PCA linear transformation is obtained by computing the eigendecomposition of the covariance matrix $K$ defined by 
\begin{equation*}
    K =\frac{1}{n-1} (X-\Bar{x})^T (X-\Bar{x})\qquad\text{where }\qquad\Bar{x}=\frac{1}{n}\sum_{i=1}^n x_i\quad.
\end{equation*}
The matrix $K$ is a symmetric matrix by definition and hence is diagonalizable by an orthogonal basis.
Therefore we can find an orthogonal matrix $P$, and a diagonal matrix $\Lambda$, such that:
$ K= P \Lambda P^T $.
Thus, the matrix $K$ defines, via $P$ and $\Lambda$, an \emph{orthonormal eigensystem}  $\left\{(\lambda_{i},v_{i})\right\}_{i=1}^{d}$, where $K v_{i}=\lambda_{i} v_{i}$, with
\begin{equation}
\label{order}
    \lambda_{1}\geq \lambda_{2}\geq \cdots \geq\lambda_{d}.
\end{equation}

To compute the orthonormal basis in $\mathbb{R}^d$ whose directions maximize the variance, PCA finds the first principal component, which is the direction in $\mathbb{R}^d$ along which projections have the largest variance among all possible directions and then iterates.
In particular, 
\begin{equation}
\label{equations}
    v_1=\max_{v\in S^{d-1}}\sum_{x\in X} (x-\Bar{x})\cdot v\qquad \rm{and }\qquad v_i=\max_{v\in\left(S^{d-1}\perp\Span\{v_j\}_{j<i}\right)}\sum_{x\in X} (x-\Bar{x})\cdot v
\end{equation}
where $S^{d-1}\subset\mathbb{R}^d$ in the standard manner. The second principal component vector, $v_2$, is the direction that maximizes variance among all directions that are orthogonal to the first principal component, with similar properties for $v_i, 1<i\leq d$.
Nontrivially, the vectors $\{v_i\}_{i=1}^d$ defined in (\ref{equations}) are precisely the eigenvectors of the covariance matrix, $K$, we refer the reader to~\cite{shlens2014tutorial} for further details.

Each eigenvalue, $\lambda_i$, from above explains the variance associated with its paired eigenvector.
As noted in~(\ref{order}) the highest eigenvalue denotes to the direction of highest variance in the data, which is the first principal component, $v_1$.
Correspondingly, the eigenvalues obtained from the covariance matrix are often referred to as the \textit{explained variances}.
The \textit{normalized explained variance} $\overline{\lambda_i}$ is defined by:
$
    \overline{\lambda_i}=\frac{1}{\sum_{j=1}^d
\lambda_j }\lambda_i.
$

\begin{remark}
\label{rem:ha}
The value $\overline{\lambda_i}$ describes the percentage of variance in the $i^{th}$ direction.
Intuitively, the normalized explained variance provides a measure of the degree of importance for each corresponding eigenvector.
As these eigenvalues are ordered~(\ref{order}), one can also refer to the most important $k$ PCA directions, where $k\leq d$, and these are the vectors $v_1,\ldots,v_k$ corresponding to the largest $k$ eigenvalues.
These properties of the PCA orthonormal eigensystem $\left\{(\lambda_{i},v_{i})\right\}_{i=1}^{d}$ can be used to obtain a heuristic assessment of the dimensionality of distinguishing features within the original data set $X$.
This heuristic is obtained by measuring the cumulative values of the normalized explained variance $S_k:= \sum_{i}^k \overline{\lambda_i} $, where $1 \leq k \leq d$. In practice, we choose $k$ such that $S_k \geq r$ where $r$ is a reasonable percentage value for which the chosen PCA vectors still capture the original data. In this work we have chosen $r$ to be $0.95$. 
\end{remark}

\section{Preparing the Data}
\label{sec:data}
Regarding the data sets used in this paper we note the following.
Most knots exist in pairs, $(K,\mir(K))$, such as the left and right handed trefoil knots, where one knot becomes the other when all positive and negative crossings are switched, so $L_+\longleftrightarrow L_-$ everywhere, giving what is referred to as the mirror image of the knot.
Knot tabulation efforts have generally accepted that it is not necessary to enumerate both of these paired knots when listing all knots, but many invariants are sensitive to this choice.
For the Jones polynomial $J_{\mir(K)}(q)=J_K(q^{-1})$~\cite{lickorish2012introduction,jones1997polynomial}, while for the signature $\sigma(\mir(K))=-\sigma(K)$~\cite{lickorish2012introduction}, and for the Rasmussen s-invariant $s(\mir(K))=-s(K)$ as well~\cite{rasmussen2010khovanov}.

In light of this symmetry, to reduce memory overhead and computation time, as well as to enhance the clarity of the associated data visualizations, we have chosen to include just one of either $K$ or $\mir(K)$ in our data set.
We first chose the embedding where the signature was positive. If the signature was zero, we then chose $K$ or $\mir(K)$, to ensure the Rasmussen s-invariant was positive. When both the signature and s-invariant are zero, or one was unknown, we chose the embedding of $K$ or $\mir(K)$ for which the most extreme degree was positive (e.g. by choosing $\mir(6_1)$ over $6_1$ as $J_{\mir(6_1)}(q)=q^4-q^3+q^2-2 q+2-\frac{1}{q}+\frac{1}{q^2}$ has extreme degrees of 4 and -2) as in Table~\ref{tab:jonesVectors}.
Once a choice between $K$ and $\mir(K)$ was made we constructed each point cloud by the following general method.

Given a finite family of knots $\mathcal{F}$, and a single variable knot polynomial invariant $I$, we want to construct a point cloud $\mathcal{C}_{\mathcal{F}}^I\subset  \mathbb{R}^n$. The procedure we describe here works for any finite set of knots and any single variable polynomial knot invariant $I$. The steps for creating this point cloud are as follows: 
\begin{enumerate}
    \item For each $K \in \mathcal{F}$ we compute the polynomial invariant $I(K)$ or $I(\mir(K))$ as in the second column of Table~\ref{tab:jonesVectors}.
    \item Convert each polynomial $I(K)$ to a coefficient vector $P(I(K))$ - or $P(I(\mir(K)))$.
    \item The set of coefficient vectors $\{P(I(K)) \mid K\in\mathcal{F}\}$ are aligned by padding each vector with zeroes to ensure that the coefficient of $q^0$ is in a consistent position and all vectors are of the same length as in the right column of Table~\ref{tab:jonesVectors} and elaborated upon below. Each padded vector is denoted $P_{\mathcal{F}}(I(K))$.
\end{enumerate}

An example of this method is given for a small set of knots in Table~\ref{tab:jonesVectors}. There we present the choice of embedding for each knot, corresponding Jones polynomial and the resulting vector in the point cloud of the family of knots up to 6 crossings.

\begin{table}[!ht]
\begin{tabular}{c|l||r@{~}r@{~}r@{~}r@{~}r@{~}r@{~}r@{~}r@{~}r@{~}r@{~}r}
$K$& $J(K)$&&&&$q^0$\\
\hline
$0_1$&1
& (~0,&0,&0,&1,&0,&0,&0,&0,&0,&0,&0)\\
$\mir(3_1)$& $q+q^3-q^4$
& (~0,&0,&0,&0,&1,&0,&1,&-1,&0,&0,&0)\\
$4_1$&  $q^{-2}-q^{-1}+1-q+q^2$
&(~0,&1,&-1,&1,&-1,&1,&0,&0,&0,&0,&0)\\
$\mir(5_1)$& $q^2+q^4-q^5+q^6-q^7$
&(~0,&0,&0,&0,&0,&1,&0,&1,&-1,&1,&-1)\\
$\mir(5_2)$& $q-q^2+2q^3-q^4+q^5-q^6$
&(~0,&0,&0,&0,&1,&-1,&2,&-1,&1,&-1,&0)\\
$\mir(6_1)$&$q^{-2}-q^{-1}+2-2q+q^2-q^3+q^4$
& (~0,&1,&-1,&2,&-2,&1,&-1,&1,&0,&0,&0)\\
$\mir(6_2)$& $q^{-1}-1+2q-2q^2+2q^3-2q^4+q^5$
&(~0,&0,&1,&-1,&2,&-2,&2,&-2,&1,&0,&0)\\
$6_3$&$-q^{-3}+2q^{-2}-2q^{-1}+3-2q+2q^2-q^{3}$
&(-1,&2,&-2,&3,&-2,&2,&-1,&0,&0,&0,&0)
\end{tabular}
\caption{Jones polynomials with positive extreme degree and their padded coefficient vectors aligned at $q^0$.}
\label{tab:jonesVectors}
\end{table}

Observe that each vector $P_{\mathcal{F}}(I(K))$ does not solely depend on the knot, $K$, but rather also depends on the family $\mathcal{F}$.
Indeed, the coefficient vectors, $P(I(K))$ for various knots frequently are of differing lengths and belong to different Euclidean spaces.
Constructing $\mathcal{C}_{\mathcal{F}}^I$ depends explicitly on the family $\mathcal{F}$, since we padded the shorter vectors in this set to match the longest ones.
Even when two polynomials in a family have the same length coefficient vector as with $P(I(\mir(5_1)))$ and $P(I(\mir(5_2)))$ in Table~\ref{tab:jonesVectors} on the left, alignment frequently pads the vectors differently.
In order to obtain the point cloud data $\mathcal{C}_{\mathcal{F}}^J$, where $\mathcal{F}$ is all knots with at most 6 crossings, we align the vectors as shown in the right column of Table~\ref{tab:jonesVectors}, padding with zeros as necessary.

For generating the data here
we used the KnotTheory package~\cite{KA} to compute the signature and the Jones polynomial for all knots up to 17 crossings and received a collection of the s-invariants for all knots up to 15 crossings from Alex Shumakovitch~\cite{Priv_Comm_Person_PCA_Model}. 

For each invariant, we then defined coefficient vectors for each polynomial, padding the vectors with zeroes as needed to align the position of $q^0$ as done to go from left to right in Table~\ref{tab:jonesVectors}.

\section{The Filtration Method}
\label{filter}

In data analysis one commonly has to draw conclusions based on incomplete data sets.
Therefore, one hopes to find those properties that do not evolve, but rather persist unchanged throughout the data set, suggesting that the conclusions are fundamental to all of the data rather than as a property of whichever special subset is being considered.
To address this issue, we have applied manifold learning to filtrations of our data, that is, nested sequences of point clouds indexed with respect to some increasing parameter.
This allows us to both detect essential, conjecturally constant, features, while also detecting those that meaningfully evolve with respect to different parameters.

\begin{definition}
\label{nested}
A \emph{
filtration} of a set $\mathcal{F}$ is a finite sequence $\{\mathcal{F}_i \}_{i=1}^n$ of nested  sets such that: 
   $ \mathcal{F}_1 \subset \mathcal{F}_2 \subset \cdots \subset  \mathcal{F}_n=\mathcal{F}.$ 
\end{definition}

For our purposes, let each set $\mathcal{F}_i $ be a family of knots, and $I$ be a single variable polynomial knot invariant.
The nested sequence in Definition \ref{nested} induces a nested sequence on a corresponding filtration of point clouds, denoted: 
   \begin{equation}
\label{order point cloud}
     \mathcal{C}^{I}_{\mathcal{F}_1} \subset \mathcal{C}^{I}_{\mathcal{F}_2} \subset \cdots \subset \mathcal{C}^{I}_{\mathcal{F}_n}.
\end{equation}
Now even though $\mathcal{F}_j\subset \mathcal{F}_{j+1},$ the corresponding vectors $P_{\mathcal{F}_j}(I(K))$ and $P_{\mathcal{F}_{j+1}}(I(K))$ often belong to different Euclidean spaces, but there is always a natural mapping that sends a point in $P_{\mathcal{F}_j}(I(K))$ to the corresponding point in $P_{\mathcal{F}_{j+1}}(I(K))$.
Namely the two vectors can be aligned on the position of $q^0$, padding the necessary zeros in $P_{\mathcal{F}_j}(I(K))$ so it is embedded in the same Euclidean space as $P_{\mathcal{F}_{j+1}}(I(K))$. Using this mapping we can meaningfully embed a point cloud $\mathcal{C}^{I}_{\mathcal{F}_j}$ inside $\mathcal{C}^{I}_{\mathcal{F}_{j+1}}$ whenever $\mathcal{F}_j \subset \mathcal{F}_{j+1} $.

Studying the eigensystems generated by PCA on a nested sequence of point clouds provides insight on how their principal component vectors and corresponding normalized explained variances evolve as the size of the point cloud grows.
This can provide more information about the distribution from which this point cloud is drawn.
For instance, by considering the relative sizes between consecutive normalized explained variances from the filtration we get information on the shape of the point cloud.

There are over 50 distinct well-known invariants used to distinguish knots in tabulations~\cite{livingston}.
Most of these naturally define ordered families of knots.
One such natural filtration of families of knots 
is the one obtained by considering all knots, up to ambient isotopy, with minimal number of crossings less than or equal to a particular value $k$.
Let $\mathcal{F}_{\lceil k\rceil}$ denote the family of all different knots with crossing number less than or equal to $k$.
We wish to study the point clouds obtained by the following filtration
\begin{equation}
\label{order family2}
    \mathcal{F}_{\lceil3\rceil} \subset \mathcal{F}_{\lceil4\rceil} \subset \cdots \subset \mathcal{F}_{\lceil j\rceil}  \subset \cdots \subset \mathcal{F}_{\lceil k \rceil}. 
\end{equation}
When paired with a single variable polynomial knot invariant $I$, the filtration (\ref{order family2}) induces a filtration of point clouds as in (\ref{order point cloud}). We will refer to the point cloud filtration induced by (\ref{order family2}) as a \emph{crossing filtration}.    
To specify the PCA eigensystem obtained at each point cloud in a crossing filtration we will associate to the point cloud $\mathcal{C}^{I}_{\mathcal{F}_{\lceil j\rceil}}$ the PCA eigensystem $\left\{ \overline{\lambda_i}(j),v_i(j)\right\}$.

The second natural filtration on point clouds associated with polynomial knot invariants is induced by the norm.
We default to the $l_2-$norm due to its ease of use and scalability within the scikit-learn package, but comparisons against other norms did not produce significantly different results.
Given $\mathcal{F}$, a finite family of knots with $I$, a single variable polynomial knot invariant, then for any $r\in \mathbb{R}^+ $, we define
\begin{equation*}
    \mathcal{C}^{I}_{\mathcal{F}}( r ) :=\{ p\in\mathcal{C}^{I}_{\mathcal{F}} \mid \norm{p}\leq r\}.
\end{equation*}
For a finite sequence of positive real numbers $r_1 \leq \cdots \leq r_{max}$, we define the filtration:
\begin{equation}
 \label{norm filter}   
\mathcal{C}^{I}_{\mathcal{F}_{\lceil k\rceil}}( r_1)\subset\mathcal{C}^{I}_{\mathcal{F}_{\lceil k\rceil}}( r_2)\subset\mathcal{C}^{I}_{\mathcal{F}_{\lceil k\rceil}}( r_3)\subset\cdots\subset\mathcal{C}^{I}_{\mathcal{F}_{\lceil k\rceil}}( r_{max})=\mathcal{C}^{I}_{\mathcal{F}_{\lceil k\rceil}},
\end{equation}
where $r_{max} = \mathrm{max}_{p\in\mathcal{C}^{I}_{\mathcal{F}_{\lceil k\rceil}}}\left(\norm{p}\right).$
We call the point cloud filtration given in (\ref{norm filter}) the \textit{norm filtration} and will denote the PCA eigensystem obtained from $\mathcal{C}^{I}_{\mathcal{F}_{\lceil k\rceil}}( r_j)$ by $\left\{ \overline{\lambda_i}(r_j),v_i(r_j)  \right\}$. 

\section{Application to the Jones polynomial}
\label{sec:results}

In this section we outline the results of running PCA on point clouds obtained from the Jones polynomial.
In total we use all 9,755,329 knots of at most 17 crossings.
In Subsection~\ref{ssec:norm} we discuss how filtering by the $l_2$-norm illustrates the shape of the point cloud $\CF{J}{17}$ using the filtration discussed in Section~\ref{filter}.
While in Subsection~\ref{ssec:crossing} we examine how the crossing filtration illuminates persistent features of the point cloud.
To see if this behavior continues at higher crossing number we consider the point cloud for two special subfamilies of knots up to 2001 crossings in Section~\ref{sec:special}.

Our analysis utilizes two primary visualizations of the PCA data.
In the first, the relative importance of each PCA eigenvector can be visualized by plotting the normalized explained variance as a function of the ordered indices of the eigensystem as discussed in Remark~\ref{rem:ha}.
Similarly, one can also plot the cumulative values of the normalized explained variance $S_k$ as a function of $k$. 
Our second visualization follows the trajectory of a sequence of PCA eigensystems $\left\{ \overline{\lambda_i},v_i \right\}_j$ across each step of the filtration by comparing the values of $\overline{\lambda_i}$ and the directions of the principal components $v_i$.
Ideally, all the principal component vectors would overlap across the filtration, we measure their deviation by using the classical dot product between them:
\begin{align}
v_i(r_{j+1})\cdot v_{i}(r_j)&=\cos{\theta_{i,j}}\norm{v_i(r_{j+1})}\norm{v_{i}(r_j)}\qquad\text{in the norm filtration;}\label{eqn:normAngle}\\
v_i(j+1)\cdot v_{i}(j)&=\cos{\theta_{i,j}}\norm{v_i(j+1)}\norm{v_{i}(j)}\qquad\text{in the crossing filtration.}
\end{align}
Here $\theta_{i,j}$ is the angle between the $i^{th}$ principal component in the $j$ and $j+1^{st}$ eigensystem.

\subsection{Structure from the Norm Filtration}
\label{ssec:norm}

\begin{figure}[!ht]
  \centering
   {\includegraphics[scale=0.49]{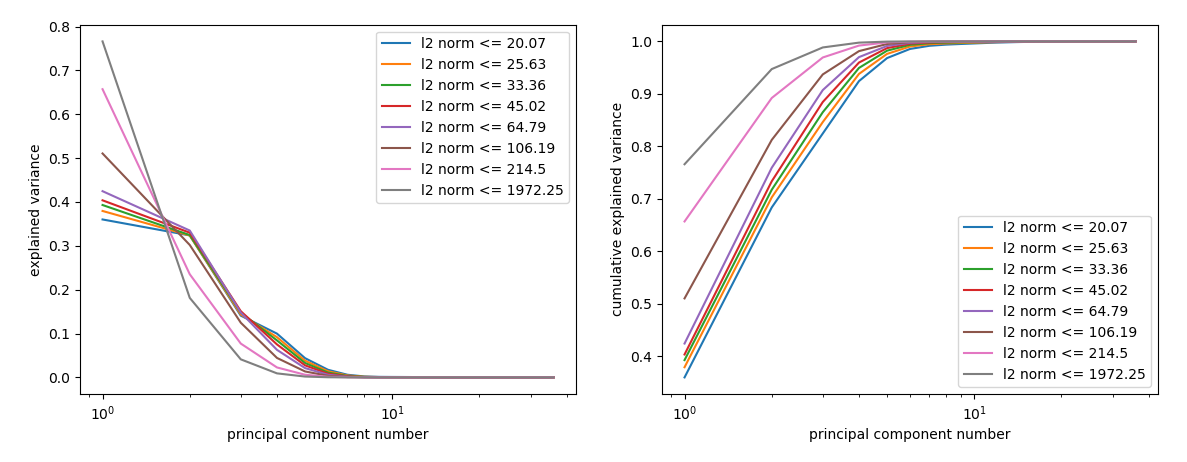}
    \caption{
    The normalized explained variance and cumulative normalized explained variance of each principal component in the norm filtration on the Jones polynomial, plotted on a log scale to highlight the  behavior of the most significant components.
    } \label{Jones 17}
  }
\end{figure}

Having prepped and filtered the data as in Sections~\ref{sec:data} and~\ref{filter}, we set up the norm filtration of (\ref{norm filter}) using the radii $\{r_7,\ldots,r_i,\ldots,r_0\}$, denoted in Figure~\ref{Jones 17}.
Each radius is chosen to restrict to the central  $\frac{1}{2^i}$ of the point cloud, doubling the number of points considered with each iteration.

Figure \ref{Jones 17} illustrates the first type of visualization, where we consider the plotted values of the normalized explained variance.
The left hand set of curves show the normalized explained variance of each principal component, ordered as in~(\ref{order}), while the right hand graph shows the cumulative normalized explained variance.
In both graphs the PCA calculation is done on the filtration 
$\mathcal{C}^{J}_{\mathcal{F}_{\lceil 17\rceil}}( r_7)\subset\mathcal{C}^{J}_{\mathcal{F}_{\lceil 17\rceil}}(r_6)\subset\cdots\subset\mathcal{C}^{J}_{\mathcal{F}_{\lceil 17\rceil}}(r_0)$ with each family denoted by a distinct color.
The exponential division of the point cloud is both for eventual contrast with the exponential growth of the crossing filtration and to ensure that at each step the amount of new data is equal to the amount of preexisting data.

Two trends quickly appear in the Figure \ref{Jones 17}.
The first principal component becomes more significant as the bounding radius of the point cloud increases, while subsequent components decrease in prominence.
This increase in the first principal component exceeds the decrease in subsequent components and as a result the cumulative normalized explained variance, $S_k$, increases for each $k$.
Following the bound set out in Remark~\ref{rem:ha} we see that for $r_7,r_6,$ and $r_5$ $S_5\geq0.968$, while for $r_4,r_3,$ and $r_2$, $S_4\geq0.959$, and then $r_1$ and $r_0$ have $S_3\geq0.969$ and $S_3\geq0.988$ respectively, which suggests that $\mathcal{C}^{J}_{\mathcal{F}_{\lceil 17\rceil}}$ approximates a 3-dimensional manifold.

\begin{figure}[!ht]
  \centering
   {\includegraphics[scale=0.46]{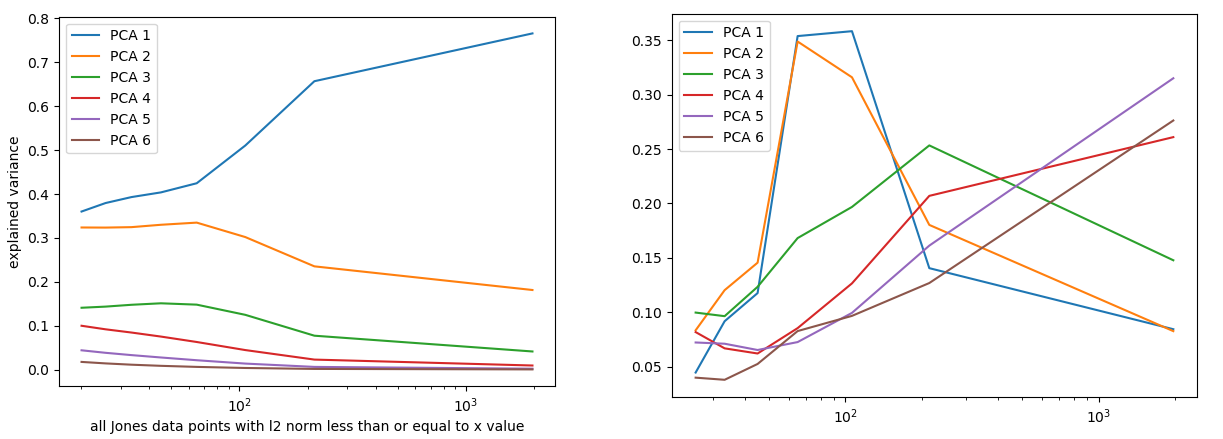}
    \caption{
    The left hand chart plots the normalized explained variance against the radius of the norm filtration of $\mathcal{C}^J_{\mathcal{F}_{\lceil17\rceil}}$ for the first $6$ principal components on a log scale in accordance with our doubling sample size. Here $\overline{\lambda_1}$ grows as the size of this knot family increases, while the $\overline{\lambda_i}$ for $i\geq 2$ principally decrease. The figure on the right provides insight on any trends toward stability in the angle $\theta_{i,j}$. In this notation, the $x$ axis represents the log of the bounding norm of $r_j$, while the $y$ axis represents the angle $\theta_{i,j}$ (measured in radians) and the index $i$ is depicted using different colors.
    }
    \label{Jones 17 2}
  }
\end{figure}

These trends are affirmed by the second type of visualization in Figure~\ref{Jones 17 2}, where the trajectory of the first six components of the normalized PCA eigensystem are followed across each step of the filtration.
The left graph illustrates the gradual growth of the first principal component at the expense of the remaining components as the radius of our point cloud  $\mathcal{C}^J_{\mathcal{F}_{\lceil17\rceil}}(r_i)$ increases.
Quantitatively, we solely note that the maximum relative spread for any of the 3 significant PCA components is $\sim91\%$.

The quality of a filtration is not only reflected in the stability of the $\overline{\lambda_j}$'s, we can also measure the alignment of sequential eigenvectors as in (\ref{eqn:normAngle}).
On the right hand side of Figure~\ref{Jones 17 2} we plot these angles for principal components $1\leq j\leq6$ across radii $15< r< 2000$.
The principal components stabilize as the filtration radius increases, but two details stand out.
First, from $r_4$ to $r_2$ the variation between important eigenvectors reduces to a stable point.
Secondly, the angles between secondary sequential components begins to stabilize from $r_2$ to $r_0$.
Furthermore, the significance of the principal component does not appear to correlate directly with the relative degree of stability across the filtration.

Figure~\ref{Jones 17 2} suggests that in the center of the point cloud, the data spreads fairly evenly in 2 directions, before changing direction between radii of 50-200 and pronouncedly extending out in a single new direction.
The cutoffs in the data based on the doubling radii of the point cloud suggests that the data is disproportionately densely packed towards the center of the distribution and is sparse towards the extremes.
We consider the shape of the data further when discussing Figure~\ref{alternating VS all knots}.

\subsection{Persistent Properties in the Crossing Filtration}
\label{ssec:crossing}
\begin{figure}[!ht]
  \centering
   {\includegraphics[scale=0.49]{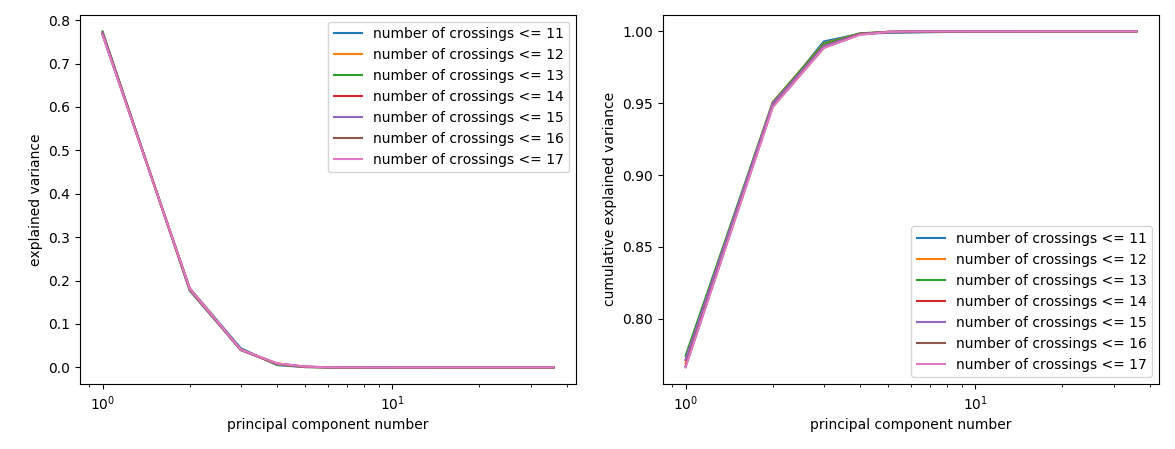}
    \caption{
    The normalized explained variance and cumulative normalized explained variance of each principal component in the crossing filtration on the Jones polynomial plotted on a log scale to emphasize the behaviour of the most important components.
    }
    \label{Jones 17 number of crossings}
  }
\end{figure}

Following the same steps we now consider the crossing number filtration following from~(\ref{order family2}). 
This filtration presents distinct features from the norm filtration.
The number of knots in each step of the filtration increases exponentially~\cite{ernst1987growth}, so to ensure a sufficient number of data points in the smallest filtration, we only consider the cases $\mathcal{C}^J_{\mathcal{F}_{\lceil11\rceil}}\subseteq\mathcal{C}^J_{\mathcal{F}_{\lceil12\rceil}}\subseteq\mathcal{C}^J_{\mathcal{F}_{\lceil13\rceil}}\subseteq\mathcal{C}^J_{\mathcal{F}_{\lceil14\rceil}}\subseteq\mathcal{C}^J_{\mathcal{F}_{\lceil15\rceil}}\subseteq\mathcal{C}^J_{\mathcal{F}_{\lceil16\rceil}}\subseteq\mathcal{C}^J_{\mathcal{F}_{\lceil17\rceil}}$.
The visualization of Figure~\ref{Jones 17 number of crossings} when compared to Figure~\ref{Jones 17} presents a strong contrast.

\begin{figure}[!ht]
  \centering
   {\includegraphics[scale=0.49]{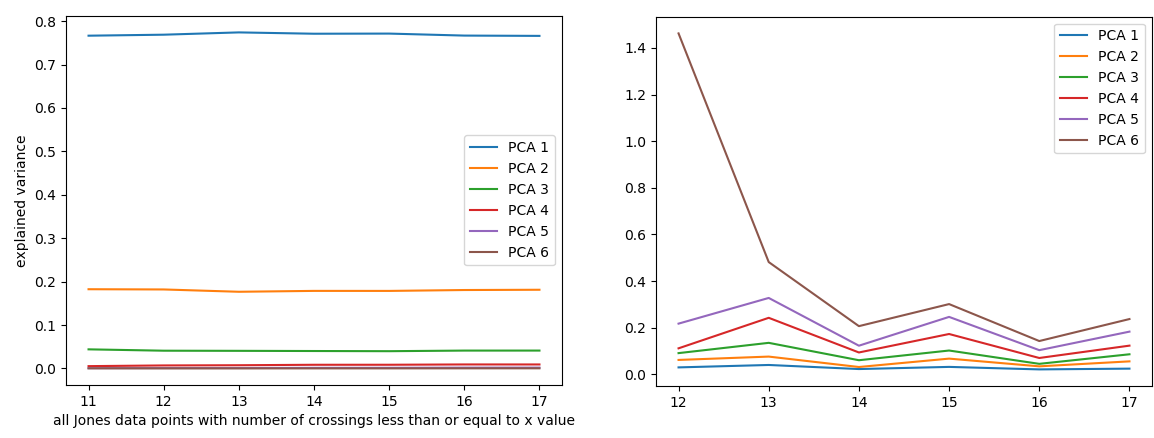}
    \caption{
    The left hand chart shows the explained variance plotted against the crossing filtration on $\CF{J}{17}$.
    It is remarkably level across the filtration.
    The right hand chart shows the stability of the principal components.
    Notably, the more significant the component, the more stable it is.} \label{Jones 17 number of crossings second}
  }
\end{figure}

In Figure~\ref{Jones 17 number of crossings} the normalized explained variances and cumulative normalized explained variances are essentially indistinguishable.
Following the bound set out in Remark~\ref{rem:ha}, we find that $0.992\geq S_3\geq 0.988$ for every family in the filtration and that $S_2<0.95$ except for the 12 and 13 crossing families, where $S_2\Dot{=}0.9507$.
We also consider the crossing number filtration analogue of Figure~\ref{Jones 17 2} in Figure~\ref{Jones 17 number of crossings second}.

Figure~\ref{Jones 17 number of crossings second} stands in marked contrast to its analogue for the norm filtration.
The normalized explained variance is remarkably stable for the $3$ significant components with a maximal relative spread of $\sim 3.5\%$, a significant improvement in consistency. The principal components also behave differently when measuring the angle between $v_i\in\CF{J}{k}$ and $v_i\in\CF{J}{k+1}$ compared to the norm filtered case.
There may be more total variation across the filtration, but that variation is more orderly, with $\theta_{i,k}<\theta_{i+1,k}$ for all $k$ and all significant $i$, where $\overline{\lambda_i}>0.00001$.
It is not surprising that these filtrations have some amount of variation as the minimal dimension of $\mathbb{R}^n$ in which $\CF{J}{k}$ can be embedded strictly increases with $k$.
Of further interest is a mild periodicity in the variation of $\theta_i$ across k, suggesting that the change in distribution of knots in $\CF{J}{k}$ depends on the parity of $k$.

\begin{figure}[!ht]
  \centering
   {\includegraphics[scale=0.205]
   {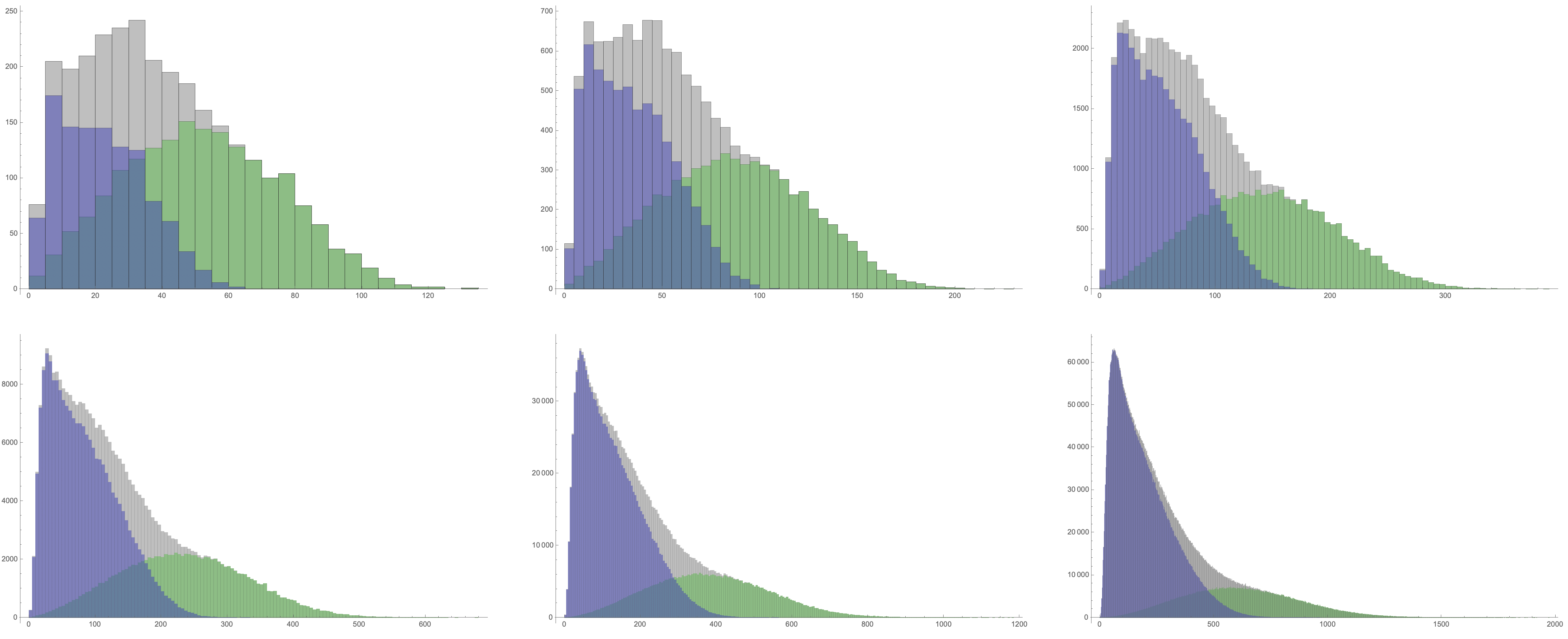}
    \caption{The distribution of the $l_2$-norms (total count vs. norm) for the Alternating (Green), NonAlternating (Blue), and Combined (Grey) knots up to 12, 13, 14, 15, 16, and then 17 crossings when taken left to right, top to bottom.
    } \label{alternating VS all knots}
  }
\end{figure}

The disparity in stabilization behavior between the norm filtration and the crossing number filtration begs the question of whether there is something special about either one.
Looking at the distribution of norms for $\CF{J}{17}$ as in the lower right of Figure~\ref{alternating VS all knots} it becomes immediately apparent that the norm filtration suffers from some structural deficiencies.
Namely, the subfamily of nonalternating knots has a skewed distribution towards lower norms, while the $l_2$-norms of alternating knots favor a broader distribution.
We observed in talking about Figure~\ref{Jones 17 2} that the angles between principal components experienced an inflection and rapid stabilization in the three most important PCA components between $r_4\doteq45$ and $r_1\doteq215$.
Additional experimentation has shown that these characteristics stabilize almost completely for $r\geq 1000$.
In Figure~\ref{alternating VS all knots} we see that, for every family in the crossing filtration, the nonalternating knots contribute an insignificant fraction of new data points to the PCA calculation by the point where the alternating knot distribution peaks.
Similarly, the tail of this distribution continues stretching as the crossing number increases, so in each case only a small number of data points are added to the point cloud after an $r<<r_{\max}$ so it is of little surprise that the principal components mostly stabilize after a given point.

This dependence on the norm distributions of the alternating and nonalternating knots suggests that we should also consider these knot classes by themselves.
Let $\CFn{I}{k}$ denote the point cloud of \emph{nonalternating knots} of at most $k$ crossings built using the single polynomial invariant $I$, and  $\CFa{I}{k}$ will denote the analogous point cloud of \emph{alternating knots}.

\begin{figure}[!ht]
  \centering
   {\includegraphics[scale=0.5]{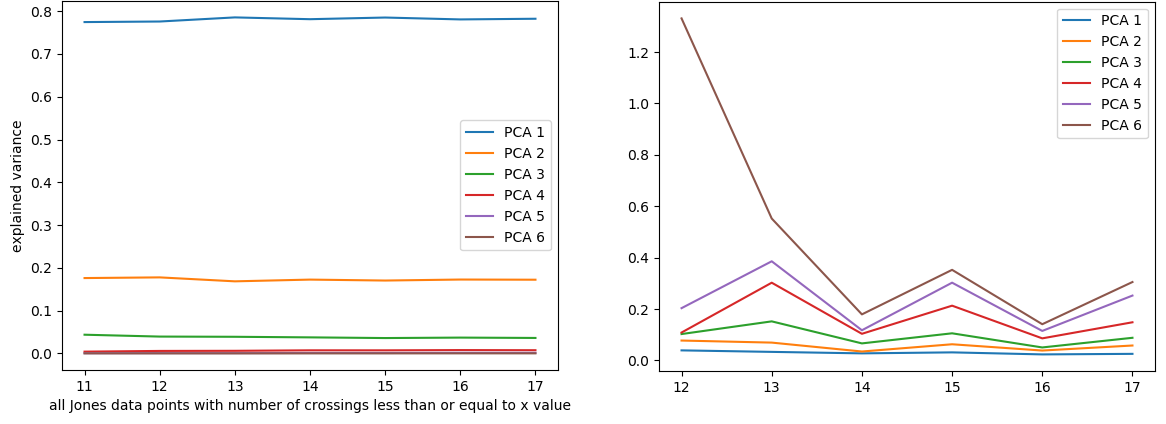}
    \caption{
    The left figure plots the normalized explained variance against the radius of the $l_2$-norm filtration of $\CFa{J}{17}$. The figure on the right shows how the PCA bases obtained from the filtration stabilize as we 
    increase the radius.
    (Note that the larger the contribution the less the deviation) } \label{alternating knots second}
  }
\end{figure}

In Figure~\ref{alternating knots second} we first consider the persistence of the PCA eigensystem features under the crossing number filtration of  $\CFa{J}{11}\subset\CFa{J}{12}\subset\cdots\subset\CFa{J}{17}$ on alternating knots.
The normalized explained variances on the left of Figure~\ref{alternating knots second}, suggests they follow the same general pattern as expressed in Figure~\ref{Jones 17 2}, but with a value of $\overline{\lambda_1}\doteq0.782$ and relative spread of $\sim 1.4\%$.

\begin{figure}[!ht]
  \centering
   {\includegraphics[scale=0.43]{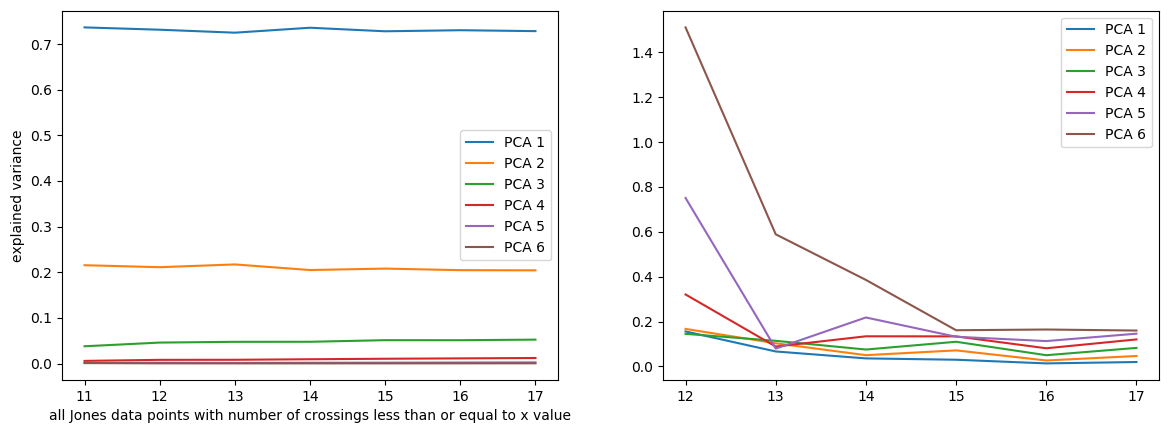}
    \caption{
    The left figure plots the normalized explained variance of $\CFn{J}{17}$. The figure on the right shows the how the PCA bases obtained from the filtration change with crossing number.
    Note that here the trends are consistent at high crossing number, but not at low crossing number where a dearth of examples likely lead to noise.}
    \label{nonalternating knots second}
  }
\end{figure}

Next we consider the persistence of the PCA eigensystem features under 
the crossing number filtration of $\CFn{J}{11}\subset\CFn{J}{12}\subset\cdots\subset\CFn{J}{17}$ on nonalternating knots,  as in  Figure~\ref{nonalternating knots second}.
Like the crossing filtration on alternating knots, the PCA eigensystem values for the crossing filtration on nonalternating knots are stable, but with $\overline{\lambda_1}\doteq0.728$ and relative spread of $\sim 1.6\%$.
It is worth noting that the normalized explained variances of the alternating knots and nonalternating knots settle at different values, but their combination, at steadily diverging weights, as illustrated by the relative proportions in Figure~\ref{alternating VS all knots}, still remains not just consistent as noted by Figure~\ref{Jones 17 number of crossings second}, but has even less relative spread with $\overline{\lambda_1}\doteq0.766$ and relative spread of $\sim 1.0\%$.
Even considering our final normalized eigenvalue deemed significant by Remark~\ref{rem:ha}, $\overline{\lambda_3}$, we find that the relative spread in all knots is $\sim 3.5\%$, while the alternating and nonalternating knots have relative spreads of  $\sim 9.2\%$ and $\sim 13\%$ respectively, when considering crossing filtrations for $12\leq k\leq 17$ to ensure at least 1000 knots in every filtration.
The possible implications of these observations bear further investigation.

\section{Examining Jones Structure at Higher Crossing Number}
\label{sec:special}

\begin{figure}[!ht]
  \centering
  { \includegraphics[height=1.22in]{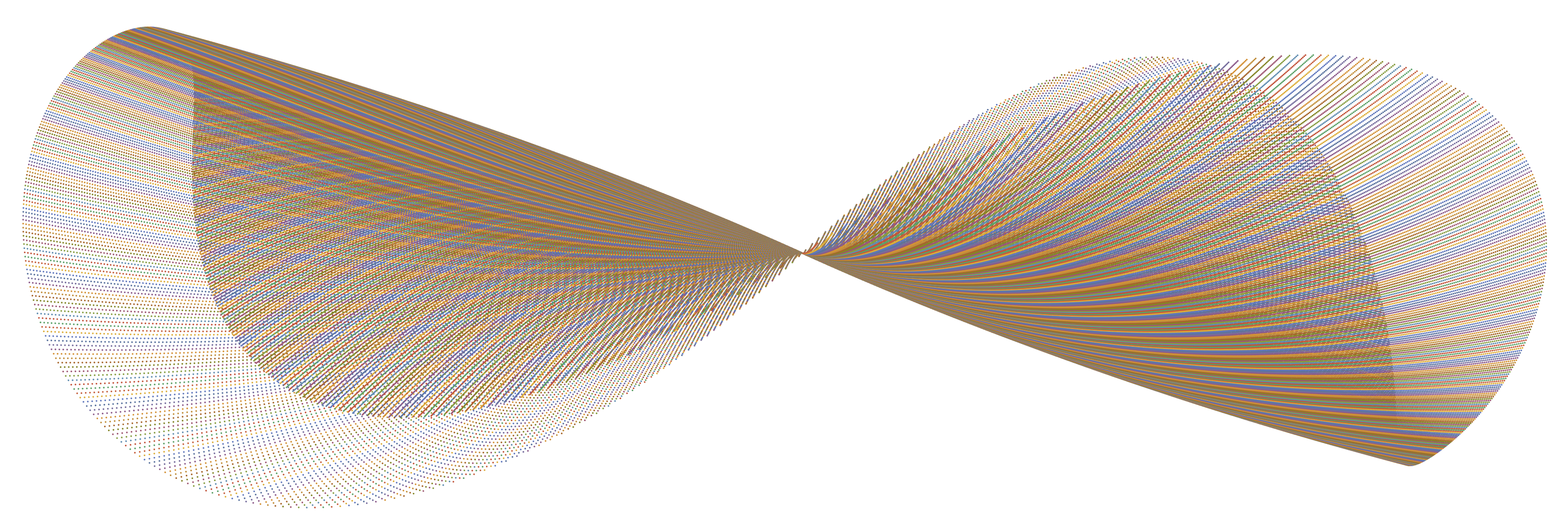} \hspace{0.5cm}\includegraphics[height=1.22in]{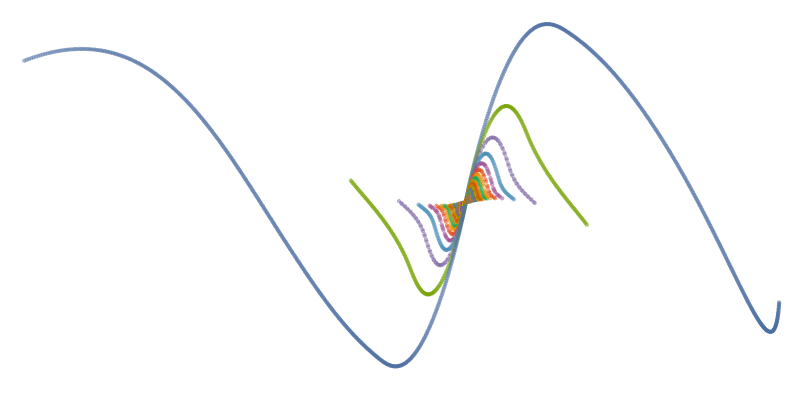}
    \caption{PCA projection into three dimensions of the positive double twist link knots up to 2001 crossings (Left) and the torus knots up to 2000 crossings (Right).}
  \label{fig:subfamilyPCA}}
\end{figure}

To provide insight into what happens for higher crossing numbers we looked at two subfamilies of knots whose Jones polynomials were easily computed at higher crossing number.
We consider the torus knots up to 2000 crossings and the positive double twist links of up to 2001 crossings.
Let $\CFt{J}{k}$ denote the point clouds of \emph{torus} knots up to $k$ crossings, and $\CFp{J}{k}$ the point clouds of single strand \emph{positive double twist link} knots up to $k$ crossings, which were calculated using~\cite{KA} and~\cite{elhamdadi2018twist} respectively.
The three dimensional PCA projections of $\CFp{j}{2001}$ and $\CFt{J}{2000}$ are presented in Figure~\ref{fig:subfamilyPCA} and suggest interesting structures exist.
Yet our results are inconsistent with those of Section~\ref{sec:results} and reveal more about the challenges of using manifold learning than they do specifically about the dimensions of $\CF{J}{k}$.

\begin{figure*}[!ht]
  \centering
   {\includegraphics[scale=0.4]{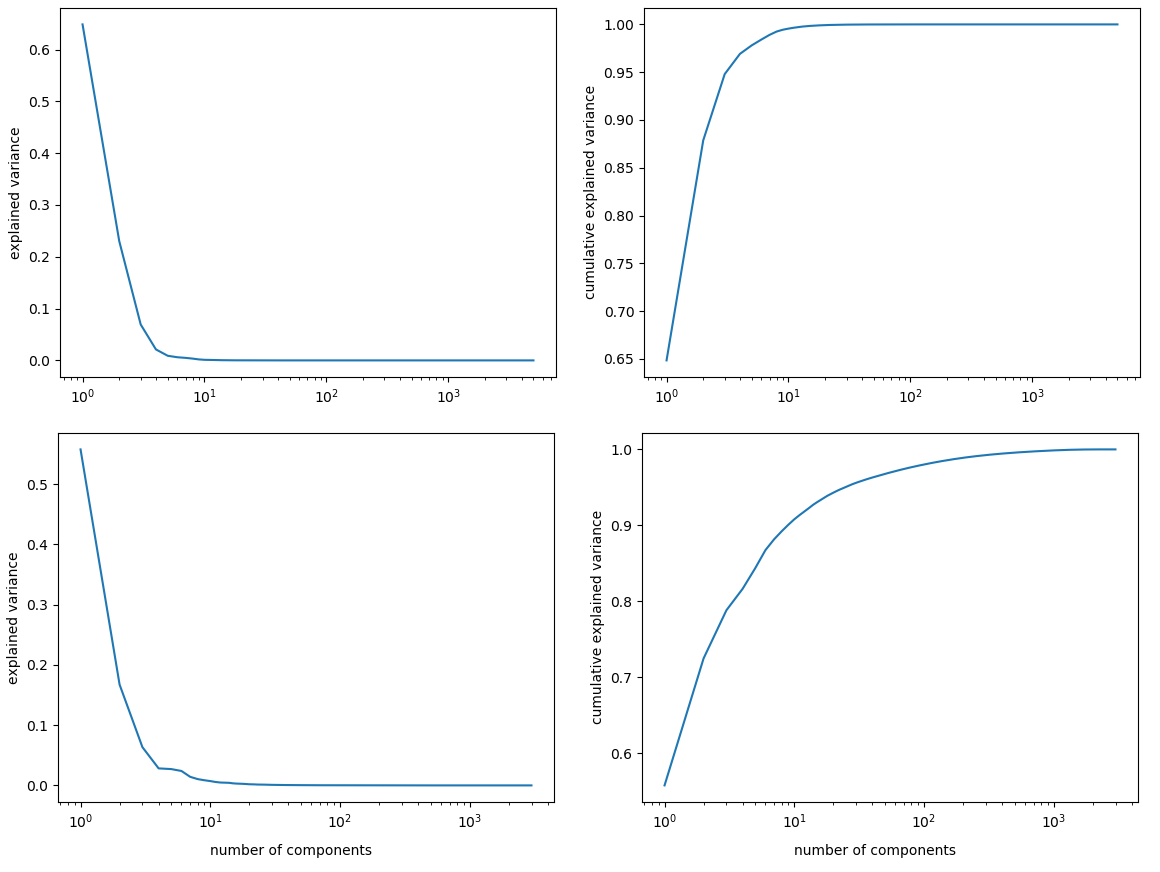}
    \caption{
    Top row: PCA on double twist link knots up to $2001$ crossings.
    Top left: $\overline{\lambda_i}$ for component $i$.
    Top right: Cumulative explained variance up to component $i$.
    Bottom row: PCA on torus knots up to $2000$ crossings.
    Bottom left: $\overline{\lambda_i}$ for component $i$.
    Bottom right: Cumulative explained variance up to component $i$.}
    \label{fig:t2p}
  }
\end{figure*}

Studying the PCA eigensystems of $\CFp{J}{2001}\subset\mathbb{R}^{5003}$ using the top row of Figure~\ref{fig:t2p} it is easy to see that $\CFp{J}{2001}$ should not be considered a 5003d manifold.
In fact, $S_4>0.969$ and $S_3\doteq0.948$, which by our heuristic suggests that $\CFp{J}{2001}$ approximates a 4 dimensional manifold.
This suggests that this submanifold of $\CF{J}{2001}$ approximates a higher dimensional manifold than we measured for $\CF{J}{17}$ and that the apparent stability in $\overline{\lambda_i}$'s seen in Figure~\ref{Jones 17 number of crossings second} might slowly evolve as crossing number increases.

We investigated this phenomenon further for torus knots.
A very different picture emerges from the analysis of the bottom row of $\CFt{J}{2000}$ in Figure~\ref{fig:t2p}.
While the left visualization is broadly similar to the results for $\CFp{J}{2001}$, the right chart displays a significant difference.
Here the cumulative normalized explained variance approaches 1 much more slowly with $S_{25}>0.95$ and taking even longer to reach a stricter restriction used by some of $S_{224}>0.99$.

Two details about $\CFt{J}{2000}$ stand out in contrast to $\CFp{J}{2001}$.
First, $\CFt{J}{2000}$ contains a mere 4501 data points unlike the over 500,000 in $\CFp{J}{2001}$.
Second, while $\CFp{J}{2001}$ lives in a 5003 dimensional space, $\CFt{J}{2000}$ lives in an 2998 dimensional space.
It is apparent that these two point clouds are not directly comparable even though they both are contained in $\CF{J}{2001}$.
This suggests that the approximate dimension of a point cloud is dependent on how it is sampled especially for nonrandom samples.
Furthermore, a direct examination of the sparsely populated $\CFt{J}{2000}$, supports the idea that for a sample size that doesn't even double the dimensionality of the space it is embedded in it is difficult to have dimensions with $\overline{\lambda_i}<<\frac{1}{\dim{\mathcal{C}^I_{\mathcal{F}}}}.$

\section{Conclusions and future work}

Studying the features of datasets that arise in pure mathematics has distinct challenges from those one faces when working with real world data.
In this paper we have outlined how to utilize one of the most traditional dimensionality reduction techniques, Principal Component Analysis, to study point clouds of data in this context.
In particular, we introduced the notion of filtrations to analyze a nested sequence of datasets.
The method introduced here is general and applicable to other scenarios where a conclusion about an infinite dataset is required.

Having explicitly described how this technique can be used to analyze the structure of the Jones polynomial data, immediate extensions of this work are to study point clouds arising from other one variable polynomial invariants such as the Alexander polynomial and to investigate the substructures illustrated in Figure~\ref{fig:pcaJones}.
In our upcoming works we will use other big data analysis techniques in the context of data in low-dimensional topology, further outlining how they can be used to compare numerical and polynomial knot invariants.
Additional dimensionality reduction calculations using ISOMAP on the $Z_0$ polynomial data affirm results obtained using PCA.
Preliminary research indicates that persistence homology confirms the existence of the substructures in the Jones polynomial data that also reflect potential relations of the Jones polynomial and signature.

\section*{Acknowledgements}
Computation for the work described in this paper was supported by the University of Southern California's Center for High-Performance Computing (\url{hpcc.usc.edu}). RS was partially supported by the Simons Collaboration Grant 318086  and NSF DMS $1854705$.

\bibliography{refs}

\end{document}